\documentclass[11pt]{amsart}
\allowdisplaybreaks
\usepackage{amsrefs} 
\usepackage{fancyhdr}
\usepackage[titletoc]{appendix}

\usepackage{amsmath,amssymb,amsthm}
\usepackage{bm}
\usepackage[dvipsnames]{xcolor}
\usepackage{etoolbox,adjustbox}
\usepackage{abstract}

\usepackage[overload]{textcase}
\usepackage{multirow}
\usepackage{graphics,graphicx}
\usepackage{float,verbatim,array}
\usepackage{todonotes}
\usepackage{xfrac,enumitem,soul}
\usepackage{tikz,tikz-cd}
\usepackage[breaklinks=true]{hyperref}
\usepackage[top=1in, bottom=1in, left=1in, right=1in]{geometry}
\usepackage{color,xcolor}
\colorlet{BLUE}{blue}
\colorlet{RED}{red}
\colorlet{GRAY}{gray}
\colorlet{BROWN}{brown}
\definecolor{OliveGreen}{rgb}{0,0.6,0}

\counterwithout{equation}{section} 
\numberwithin{equation}{section}
\newtheorem{thm}{Theorem}[section]
\newtheorem{prp}[thm]{Proposition}
\newtheorem*{theoremB}{Theorem B}

\theoremstyle{definition}
\newtheorem{dfn}[thm]{Definition}
\newtheorem{rmk}[thm]{Remark}
\newtheorem{question}[thm]{Question}
\newtheorem{example}{Example}
\newtheorem*{questionA}{Question A}

\theoremstyle{definition}
\newtheorem*{thm*}{Theorem}
\newtheorem*{lmm*}{Lemma}
\newtheorem*{crl*}{Corollary}
\newtheorem*{MRI*}{Result 1}
\newtheorem*{MRII*}{Result 2}
\newtheorem*{PRT*}{Path representation theorem}
\newtheorem*{claim*}{Claim}

\title[Set-valued metrics and generalized Hausdorff distances]{Set-valued metrics and generalized Hausdorff distances}
\author{Earnest Akofor}
\address{\textnormal{Department of Mathematics and Computer Science, 
         Faculty of Science, University of Bamenda, 
         PO Box 39 Bambili, NW Region, Cameroon}}
\email{eakofor@gmail.com}


\begin{document}

\begingroup
\def\uppercasenonmath#1{} 
\let\MakeUppercase\relax 
\maketitle

\begin{abstract}    
\vspace{0.2cm}

\noindent Let $X$ be a metric space and $BCl(X)$ the collection of nonempty bounded closed subsets of $X$. We show that Hausdorff distance $d_H$ belongs to a specific family of real-valued distances on $BCl(X)$, each of which can be expressed as the composition $\mu\circ d_{sv}$ of a topology inducing set-valued function $d_{sv}:BCl(X)^2\rightarrow \mathcal{P}(Z)$ and a real-valued set-function $\mu:\Sigma\subset\mathcal{P}(Z)\rightarrow\mathbb{R}$. With this observation, we construct several associated classes of inter-set distances, called set-valued metrics and generalized Hausdorff distances. Our constructions are both explicit and adaptable, and the resulting distance classes are expected to cover most practical applications involving distance between sets.
\end{abstract}
\let\bforigdefault\bfdefault
\addtocontents{toc}{\let\string\bfdefault\string\mddefault}
\tableofcontents
\section{Introduction}\label{Intro}
Our discussion is mainly concerned with distances between subsets of metric spaces. For generality however, we also mention subset hyperspaces. In particular, the integral Hausdorff distances of Section \ref{IntHausSec} involve subset hyperspaces that are viewed as quotients of function spaces, as recently developed in \cite{akofor2025}.

\subsection{Preliminary notation and terminology}
Let $X=(X,\tau)$ be a (topological) space. If $A\subset X$, the closure of $A$ in $X$ is denoted by $cl_X(A)$, or by $\overline{A}$ if the underlying space $X$ is understood.

Let $\mathcal{H}(X)$ denote the class of homeomorphisms of $X$ and $\mathcal{C}$ a subclass of $\mathcal{H}(X)$. By \textbf{geometry} of $X$ (respectively, \textbf{$\mathcal{C}$-geometry} of $X$) we mean the study of one or more properties of $X$ that are invariant under specific homeomorphisms of $X$ (respectively, homeomorphisms of $X$ in $\mathcal{C}$), where the invariant properties are accordingly called \textbf{geometric properties} (respectively, \textbf{$\mathcal{C}$-geometric properties}) of $X$. By \textbf{metric geometry} we mean geometry that employs metrics. With $\mathcal{P}^\ast(X)$ denoting the set of nonempty subsets of $X$, we often want to study the simplest kinds of topologies on $\mathcal{J}\subset\mathcal{P}^\ast(X)$. Since a $T_1$ space $X$ can be seen as a subset of $\mathcal{P}^\ast(X)$ in a natural way (through the inclusion $X\hookrightarrow\mathcal{P}^\ast(X),~x\mapsto\{x\}$), we often consider those topologies on $\mathcal{J}\supset\textrm{Singlt}(X):=\big\{\{x\}:x\in X\big\}$ that can be seen as extensions of the topology of $X$. Such topologies are called \textbf{hypertopologies} of $X$, and the associated spaces $\mathcal{J}$, $\textrm{Singl}(X)\subset\mathcal{J}\subset\mathcal{P}^\ast(X)$, are called \textbf{(subset) hyperspaces} of $X$. A \textbf{closed-subset hyperspace}, \textbf{compact-subset hyperspace}, or \textbf{bounded-subset hyperspace} of $X$ is a hyperspace consisting respectively of closed subsets, compact subsets, or bounded subsets of $X$. This paper is mainly concerned with certain tools for studying metric geometry of hyperspaces.

Let $X$ be a $T_1$ space, $Cl(X)\subset\mathcal{P}^\ast(X)$ a hyperspace whose points are the nonempty \textbf{closed subsets} of $X$, $K(X)\subset Cl(X)$ the subspace consisting of all nonempty \textbf{compact subsets} of $X$, and $FS_n(X):=\{A\in K(X):1\leq|A|\leq n\}$ the subspace consisting of all nonempty \textbf{finite subsets} of $X$ of cardinality at most $n$. When $X$ is a metric space, we further let $BCl(X)\subset Cl(X)$ denote the subspace consisting of all nonempty \textbf{bounded closed subsets} of $X$. Let $Y$ be a set and, as detailed in \cite{akofor2025}, let us introduce a topology $\tau_{\pi}$ on a relevant set of functions $\mathcal{F}\subset X^Y$. Since a closed set $C\in Cl(X)$ can be seen as the closure  of the image $f(Y)$ of some function $f\in X^Y$, the image-closure assignment
\begin{align}
\label{UnordMapEq}q:(\mathcal{F},\tau_{\pi})\rightarrow (Cl_Y(X),\tau_{\pi q}),~ f\mapsto cl_X(f(Y))
\end{align}
induces a quotient topology $\tau_{\pi q}$, the \textbf{$\tau_{\pi}$-quotient topology} (\textbf{footnote}\footnote{
That is, $\tau_{\pi q}:=\sup\big\{\tau~|~q:(\mathcal{F},\tau_{\pi})\rightarrow (Cl_Y(X),\tau)~\textrm{is surjective and continuous}\big\}$ is the largest topology on $Cl_Y(X)$ with respect to which the map $q:(\mathcal{F},\tau_{\pi})\rightarrow Cl_Y(X)$ is surjective and continuous.
}), on the \textbf{$Y$-indexed closed subsets}
\begin{align}
\label{IndCloSubEq}Cl_Y(X):=q(\mathcal{F}).
\end{align}
Let $X$ be a metric space and let the set of \textbf{$Y$-indexed bounded closed subsets} of $X$ be
\begin{align}
\label{IndBCloSubEq} BCl_Y(X):=Cl_Y(X)\cap BCl(X).
\end{align}

Given a point $x\in X$, a subset $A\subset X$, and a number $\varepsilon>0$, let
\begin{align*}
&\textstyle d(x,A)=\textrm{dist}(x,A):=\inf\limits_{a\in A}d(x,a),\\
&\textstyle N_\varepsilon(A):=\{x\in X:\textrm{dist}(x,A)<\varepsilon\},~\textrm{and}~\\
&\overline{N}_\varepsilon(A)=A_\varepsilon:=\{x\in X:\textrm{dist}(x,A)\leq\varepsilon\}.
\end{align*}
If $A,B\in BCl(X)$, the \textbf{Hausdorff distance} between $A$ and $B$ is
\begin{align}
\label{HausDistDef} d_H(A,B)&\textstyle:=\max\{\sup_{a\in A}\textrm{dist}(a,B),\sup_{b\in B}\textrm{dist}(b,A)\}\\
&=\inf\left\{r>0:A\cup B\subset \overline{N}_r(A)\cap \overline{N}_r(B)\right\}\nonumber\\
&=\textstyle\sup_{x\in A\cup B}|\textrm{dist}(x,A)-\textrm{dist}(x,B)|\nonumber\\
&\textstyle\stackrel{(s)}{=}\sup_{x\in X}|\textrm{dist}(x,A)-\textrm{dist}(x,B)|,\nonumber
\end{align}
where a proof of step (s) can be found in \cite[Proposition 3.1]{ConciKumb2017}.

\subsection{Motivation and related work}
Let $X$ be a metric space. In studying geometry of $BCl(X)$ as a hyperspace, we expect that the topology on $BCl(X)$ induced by a set-valued metric (representing a family of real-valued metrics by Theorem \ref{dHausCompThm}) is more general than the topology induced on $BCl(X)$ by a real-valued metric. Based on the different ways of expressing Hausdorff distance $d_H$ (Equation (\ref{HausDistDef})), our first goal is to introduce topology inducing set-valued functions (called \textbf{set-valued metrics} in Def. \ref{SVmetdfn}). Our second goal is to construct classes of real-valued generalizations of $d_H$ (called \textbf{generalized Hausdorff distances}, alias \textbf{ghd's}, in Def. \ref{HausDistClass}) by composing the above-mentioned set-valued functions with real-valued subadditive monotone set functions (called \textbf{postmeasures} in Def. \ref{PostMeasDfn}). In particular, it turns out (Theorem \ref{dHausCompThm}(iii)) that $d_H:BCl(X)^2\rightarrow\mathbb{R}$ is the composition
\[
d_H=\mu\circ d_{sv}:BCl(X)^2\stackrel{d_{sv}}{\longrightarrow}\Sigma_Z \stackrel{\mu}{\longrightarrow}\mathbb{R},
\]
of a postmeasure $\mu:\Sigma\subset\mathcal{P}^\ast(Z)\rightarrow\mathbb{R}$ and a set-valued function $d_{sv}:BCl(X)^2\rightarrow\Sigma_Z\subset\mathcal{P}^\ast(Z)$.

As indicated in Remark \ref{AdaptRmk}, our ghd-constructions are both explicit and adaptable, and the resulting distance classes are expected to cover most practical applications involving distance between sets. We note that empirical generalizations of $d_H$ exist as more suitable alternatives to $d_H$ in specific applications (see the survey of distances on sets in \cite[plus the references therein]{ConciKumb2017}, applications of distances on sets in \cite{DubJain1994,VargBog20118,BogEtal2018}, and the survey of distances in \cite{Dezas2009}).

Due to the aforementioned factoring of $d_H$, our ghd's (as members of the \emph{postmeasure family} in Definition \ref{PostMeasFam}) unify the two principal families of distances between sets (namely, \emph{Hausdorff family} and \emph{measure theoretical family}) in the survey \cite{ConciKumb2017}. Moreover, the said factoring suggests that $d_H$ is only a special case (as well as an approximation) of more general geometric tools, namely, set-valued metrics. Real-valued metrics are of course useful on hyperspaces, because they provide seemingly indispensable intermediate steps in the metric-geometric study of hyperspaces. For example, they allow relatively obscure geometric properties of a hyperspace to be revealed through metric based embedding into more familiar spaces (see \cite{Kovalev2015SP,BorovEtal2010,BorovIbra2009,chinen2015,BorsUl1931,bott1952}). Nevertheless, exclusive use of real-valued metrics on hyperspaces can ignore hidden geometric information (e.g., that involved in every possible factoring $d_H=\mu\circ d_{sv}$ in Theorem \ref{dHausCompThm}(iii) and Remark \ref{PracSvMetRmk}).

Existing notions of set-valued metrics include various notions of a hypermetric given in \cite{alim.etal2014,DeghEtal2022,NezRad2021,NezEtal2021}, which are different from our notion of a set-valued metric. Our set-valued metrics bear some resemblance to the generalized metric space structures in \cite{CookWss2021,Flagg1997,Koppmn1998}. More on the role of Hausdorff distance in the geometry of hyperspaces can be found in \cite{IllanNadl1999,LucPasq1994,zsilin1998,beer1993,Ivanova2018,PetGua2020,LeviEtal1993} and the references therein.

\subsection{Summary and highlight of main results}
The survey \cite{ConciKumb2017} studied two separate families of distances between sets, namely, the Hausdorff family and the measure theoretical family. The \textbf{Hausdorff family} consists of distances which are slight variants of $d_H$ based on various ways of expressing $d_H$. Meanwhile, with the \emph{symmetric set difference} $\Delta(A,B):=(A\cup B)\backslash(A\cap B)=A\backslash B~\cup~B\backslash A$ (for sets $A,B\in\mathcal{P}(X)$), the \textbf{measure theoretical family} consists of distances of the form $d:\Sigma(X)^2\rightarrow\mathbb{R}$, $d(A,B):=\mu\circ \Delta(A,B)$, for measures $\mu:\Sigma(X)\subset\mathcal{P}(X)\rightarrow[0,\infty]$ on $X$.
\begin{questionA}
How can we unify these two distance families based on the general form of the measure theoretical family?
\end{questionA}

In Section \ref{SvdMetSec}, by means of Definitions \ref{PartAlgDfn}, \ref{SVmetdfn}, \ref{SvOpenSetDfn}, we observe that in a metric $d:X^2\rightarrow\mathbb{R}$ on a set $X$ (Def. \ref{MetricDfn}), we can replace the ordered abelian group $\mathbb{R}=(\mathbb{R},+,\leq)$ (Def. \ref{OrdAbelGrp}) with a more general structured poset $\Sigma_Z=(\Sigma_Z,\uplus,\preceq)$, in which $\Sigma_Z\subset\mathcal{P}^\ast(Z)$ for an arbitrary set $Z$, and still be able to define a map $d_1:X^2\rightarrow \Sigma_Z$ (called a \emph{set-valued metric} in Def. \ref{SVmetdfn}) that induces a Hausdorff topology on $X$ in like manner to $d$. In conjunction with one more concept, of a \emph{postmeasure} in Def. \ref{PostMeasDfn}, we obtain the following result.

\begin{theoremB}[Theorem \ref{dHausCompThm}]
Let $M$ be a set. (i) If $d_{sv}:M^2\rightarrow\big(\Sigma_Z,\preceq,\uplus\big)$ is a sv-metric (Def. \ref{SVmetdfn}) and $\mu:(\Sigma_Z,\preceq,\uplus)\rightarrow (\mathbb{R},\leq,+)$ a real-valued postmeasure (Def. \ref{PostMeasDfn}), then the composition $\mu\circ d_{sv}:M^2\rightarrow\mathbb{R}$ is a metric.

(ii) If $\mu$ is strictly monotone {\small (i.e., $R\prec S$ implies $\mu(R)<\mu(S)$) and $\mu^{-1}(r)\backslash\{0\}\neq\emptyset$ for all $r>0$}, then $\tau_{\mu\circ d_{sv}}\subset\tau_{d_{sv}}$ (hence $\tau_{d_{sv}}$ is a Hausdorff topology).

(iii) There exists a postmeasure $\mu$ and a sv-metric $d_{sv}$ such that
\[
d_H=\mu\circ d_{sv}.
\]
 \end{theoremB}
We then answer Question A by introducing (via Theorem \ref{dHausCompThm}) the \emph{postmeasure family} (Def. \ref{PostMeasFam}) which contains both the \emph{Hausdorff family} and \emph{measure theoretical family}.

In Section \ref{GenHausSec}, we construct classes of generalized Hausdorff distances (namely, \textbf{relational} classes in Section \ref{RelHausSec} and \textbf{integral} classes in Section \ref{IntHausSec}), each of which can be associated (an example of) a set-valued metric following Theorem \ref{dHausCompThm}. The following are two examples of such distances.
\begin{enumerate}
\item A \textbf{relational} distance from Eq. (\ref{RelHausDistEq1}): For a relation~ $~R(A, B)\subset A\times B\in BCl(X)^2$,
    $$\textstyle d^R(A,B):=\sup\{d(a,b):(a,b)\in R(A,B)\}.$$
\item An \textbf{integral} distance from Eq. (\ref{PracGHDistEq2}): With $f,g\in Cl_Y(X)$ and $c_{fg}(x,y,y')$ satisfying (\ref{DIntHausEq1}) and (\ref{DIntHausEq2}),
    $$\textstyle d_{p,\nu,\mu}\left(q(f),q(g)\right):=\left(\int_{X\times Y^2}|d(x,f(y))-d(x,g(y'))|^pc_{fg}(x,y,y')d\nu_xd\mu_{y,y'}\right)^{1/p}.$$
\end{enumerate}
Necessary conditions, or even criteria, for the triangle inequality are provided (Propositions \ref{UrbGhdTIC1}, \ref{UrbGhdTIC2}, \ref{UrbGhdTIC3}, \ref{UrbGhdTIC4}) for those explicit distances for which the triangle inequality is not built-in or not obvious.

We conclude in Section \ref{DiscQuest} with some relevant questions.

\section{\textnormal{\bf Set-valued metrics}} \label{SvdMetSec}
Throughout, \emph{resp.} stands for the word \emph{respectively}. Also, ~$\overline{\mathbb{R}}_+=[0,\infty]:=[0,\infty)\cup\{\infty\}$, as the space with topology generated by open subsets of $[0,\infty)$ and the sets $\{(r,\infty)\cup\{\infty\}:r\geq0\}$.

\subsection{Construction of set-valued metrics}
\begin{dfn}[{Metric, Distance, Pseudometric, Quasimetric, Semimetric: \cite{Dezas2009}}]\label{MetricDfn}
Let $X$ be a set. A function $d:X^2\rightarrow\overline{\mathbb{R}}_+$ is a \textbf{metric} (resp., \textbf{distance}), making $X=(X,d)$ a \textbf{metric space} (resp., a \textbf{distance space}), if for all $x,y,z\in X$,
\begin{enumerate}
\item\label{MetrAx1} $d(x,y)=d(y,x)$. ~~(resp., $d(x,y)=d(y,x)$) -- symmetry.
\item\label{MetrAx2} $d(x,y)\leq d(x,z)+d(z,y)$. ~~(resp., $d(x,y)\geq 0$) -- triangle inequality.
\item\label{MetrAx3} $d(x,y)=0$ $\iff$ $x=y$. ~~(resp., $d(x,x)=0$) -- faithfulness.
\end{enumerate}
Axioms (\ref{MetrAx1})$\&$(\ref{MetrAx2}) make $d$ a \textbf{pseudometric}, (\ref{MetrAx2})$\&$(\ref{MetrAx3}) make $d$ a \textbf{quasimetric}, and (\ref{MetrAx1})$\&$(\ref{MetrAx3}) make $d$ a \textbf{semimetric}, respectively making $X=(X,d)$ a \textbf{pseudometric space}, \textbf{quasimetric space}, and \textbf{semimetric space}.
\end{dfn}

\begin{rmk}\label{SemiMetConv}~
\begin{itemize}[leftmargin=0.6cm]
\item Given any pseudometric $d_{ps}:Z^2\rightarrow \overline{\mathbb{R}}_+$, we get a metric ~$d:[Z]^2\rightarrow \overline{\mathbb{R}}_+,~([z],[z'])\mapsto d_{ps}(z,z')$~ on the set of equivalence classes $[Z]=\big\{[z]=\{z'\in Z:d(z,z')=0\}~|~z\in Z\big\}$.

\item Given any quasimetric $d_q:Z^2\rightarrow\overline{\mathbb{R}}_+$, we get a metric {\small $d:Z^2\rightarrow \overline{\mathbb{R}}_+,~(z,z')\mapsto \max\{d_q(z,z'),d_q(z',z)\}$ (or $d_q(z,z')+d_q(z',z)$)}.

\item Given any semimetric $d_s:Z^2\rightarrow\overline{\mathbb{R}}_+$, we get a pseudometric $d_{ps}:Z^2\rightarrow \overline{\mathbb{R}}_+$ given by
\begin{align}
\label{SemiPseudEq} d_{ps}(z,z')&\textstyle:=\inf\left\{\textrm{finite sums}~\sum_{i=0}^nd_s(z_{i-1},z_i):z=z_0,z_1,...,z_n=z'\in Z\right\}\nonumber\\
&\leq d_s(z,z').
\end{align}
\end{itemize}
\end{rmk}

\begin{dfn}[{Partial algebra, Partial algebra on a set}]\label{PartAlgDfn}
A \textbf{partial algebra} is a tuple $\Sigma=(\Sigma,\preceq,\uplus)$ consisting of a nonempty poset $(\Sigma,\preceq)$ and a map $\uplus:\Sigma^2\rightarrow\Sigma$ such that the following hold:
\begin{enumerate}
\item\label{PartAlgAx1} (i) $\Sigma$ contains a minimum element (\textbf{footnote}\footnote{
    A minimum element is unique: If $E,E'$ are both minimum elements then $E\preceq E'$ and $E'\preceq E$, and so $E=E'$.
    }), i.e., an element $0=0_\Sigma:=\min(\Sigma)\in\Sigma$, called \textbf{zero}, such that $0\preceq E$, for all $E\in\Sigma$. (ii) Any nonzero elements $\varepsilon,\varepsilon'\in\Sigma$ have a nonzero \textbf{common lower bound} $clb(\varepsilon,\varepsilon')\leq\{\varepsilon,\varepsilon'\}$. (iii) For any $\varepsilon\in\Sigma$, there exists $\varepsilon'\in\Sigma$ such that $\varepsilon'\uplus\varepsilon'\preceq\varepsilon$.
\item\label{PartAlgAx3} $R\uplus S=S\uplus R$.
\item\label{PartAlgAx4} If $R\preceq S$, then $R\uplus T\preceq S\uplus T$, for all $T$.

\item\label{PartAlgAx5} If $R\prec S$ and $R'\prec S'$, then $R\uplus R'\prec S\uplus S'$.
\end{enumerate}

If there is a set $Z$ such that $\Sigma\subset\mathcal{P}^\ast(Z)$, then we say $\Sigma=\Sigma_Z$ is a \textbf{partial algebra on $Z$}.
\end{dfn}
\begin{rmk}
Axiom \ref{PartAlgAx4} above is equivalent to ``If $R\preceq S$ and $R'\preceq S'$, then $R\uplus R'\preceq S\uplus S'$.''. Indeed, if Axiom \ref{PartAlgAx4} holds, then ``$R\preceq S$ and $R'\preceq S'$ imply $R\uplus R'\preceq S\uplus R'\preceq S\uplus S'$''. Conversely, if ``$R\preceq S$ and $R'\preceq S'$ imply $R\uplus R'\preceq S\uplus S'$'', then by setting $R'=S'=T$, we obtain Axiom \ref{PartAlgAx4}.
\end{rmk}

Structures that share some properties with \emph{partial algebras} include \emph{value semigroups} in \cite{Koppmn1998} and \emph{value-quantales} in \cite{Flagg1997}, with basic topological implications highlighted in \cite{CookWss2021}.

\begin{dfn}[{Set-valued metric (sv-metric), sv-distance, sv-pseudometric, sv-quasimetric, sv-semimetric}]\label{SVmetdfn}
Let $M$ be a set and $\Sigma=(\Sigma,\preceq,\uplus)$ a partial algebra. A map $d=d^{sv}:M^2\rightarrow\Sigma$ is a \textbf{set-valued metric} or \textbf{sv-metric} (resp., \textbf{sv-distance}) on $M$, making $M=(M,d)$ a \textbf{sv-metric space} (resp., \textbf{sv-distance space}) if for all $a,b,c\in M$,
\begin{enumerate}
\item\label{SvMetrAx1} $d(a,b)=d(b,a)$. ~~(resp., $d(a,b)=d(b,a)$) -- symmetry.
\item\label{SvMetrAx2} $d(a,b)\preceq d(a,c)\uplus d(c,b)$. ~~(resp., $d(a,b)\succeq 0$) -- triangle inequality.
\item\label{SvMetrAx3} $d(a,b)=0$ $\iff$ $a=b$. ~~(resp., $d(a,a)=0$) -- faithfulness.
\end{enumerate}
Axioms (\ref{SvMetrAx1}),(\ref{SvMetrAx2}) make $d$ a \textbf{sv-pseudometric}, (\ref{SvMetrAx2}),(\ref{SvMetrAx3}) make $d$ a \textbf{sv-quasimetric}, and (\ref{SvMetrAx1}),(\ref{SvMetrAx3}) make $d$ a \textbf{sv-semimetric}, respectively making $M=(M,d)$ a \textbf{sv-pseudometric space}, \textbf{sv-quasimetric space}, and \textbf{sv-semimetric space}. We will say that $d$ is \textbf{nontrivial} if there exists no order embedding $(\Sigma,\preceq,\uplus)\hookrightarrow(\mathbb{R},\leq,+)$, where a map $f:(\Sigma,\preceq,\uplus)\rightarrow(\mathbb{R},\leq,+)$ is an \textbf{order embedding} if (i) it is injective and (ii) $\forall a,b\in\Sigma$, $a\prec b$ implies $f(a)<f(b)$.
\end{dfn}

More detailed discussions on topologies of generalized metric spaces can be found in \cite{CookWss2021} and the references therein.

\begin{dfn}[{sv-metric topology, sv-metric space}]\label{SvOpenSetDfn}
Let $M$ be a set, $\Sigma=(\Sigma,\preceq,\uplus)$ a partial algebra, and $d=d^{sv}:M^2\rightarrow\Sigma$ a sv-(quasi)metric. If $0\prec\varepsilon\in\Sigma$ and $m\in M$, define the \textbf{open $d$-neighborhood} (or \textbf{$d$-ball}) of center $m$ and radius $\varepsilon$ to be
\begin{align}
\label{SvBallEqn}B_\varepsilon(m):=\big\{m'\in M:d(m,m')\prec \varepsilon\big\}.
\end{align}
The \textbf{$d$-topology} (called a \textbf{sv-metric topology}) on $M$ (making $M=(M,d)$ a \textbf{sv-metric space}) is the topology
\[
\tau_d:=\{O\subset M:\textrm{$\forall m\in O$, some $B_\varepsilon(m)\subset O$}\},
\]which is indeed a topology on $M$ by Definition \ref{PartAlgDfn} Axiom (\ref{PartAlgAx1})(ii). ({\bf footnote}\footnote{
Showing $\tau_d$ is a topology: It is clear that $\emptyset,M\in\tau_d$ and that $\tau_d$ is closed under arbitrary unions. To show $\tau_d$ is closed under finite intersections, let $A,B\in\tau_d$. If $x\in A\cap B$, then some $B_\varepsilon(x)\subset A$ and some $B_{\varepsilon'}(x)\subset B$, and so $B_{clb(\varepsilon,\varepsilon')}(x)\subset B_\varepsilon(x)\cap B_{\varepsilon'}(x)\subset A\cap B$. This shows $A\cap B\in\tau_d$.
}). By Definition \ref{PartAlgDfn} Axiom (\ref{PartAlgAx1})(iii), the  interior of a set $A\subset M$ satisfies
\[
A^o=\{m\in A:\textrm{some $B_\varepsilon(m)\subset A$}\}.~~\textrm{({\bf footnote}\footnotemark).}
\]\footnotetext{
To check this, let $O:=\{m\in A:\textrm{some $B_\varepsilon(m)\subset A$}\}$. We need to show $O=A^o$. Let $x\in O$. Then some $B_\varepsilon(x)\subset A$. Pick $\varepsilon'$ such that $\varepsilon'\uplus\varepsilon'\preceq\varepsilon$. Let $y\in B_{\varepsilon'}(x)$ and $z\in B_{\varepsilon'}(y)$. Then
\[
d(x,z)\preceq d(x,y)\uplus d(y,z)\prec \varepsilon'\uplus\varepsilon'\preceq\varepsilon,~~\Rightarrow~~B_{\varepsilon'}(y)\subset B_\varepsilon(x)\subset A,~~\Rightarrow~~y\in O,~~\Rightarrow~~B_{\varepsilon'}(x)\subset O.
\]
Therefore $O$ is the largest open set contained in $A$, because any open set contained in $A$ is also contained in $O$.
}
The interior of a ball is nonempty since $m\in B_\varepsilon(m)^o$.
\end{dfn}

\begin{example}[{Internal sv-metric}]\label{IntSvMetEx}
In Definition \ref{SvOpenSetDfn}, let $Z$ be a set, $M\subset\mathcal{P}^\ast(Z)$, $\preceq$ the inclusion $\subset$, $\uplus$ the union $\cup$, {\footnotesize $\Sigma_Z\subset(\mathcal{P}^\ast(Z),\subset,\cup)$} a partial algebra of subsets of $Z$ (where $0:=\emptyset$), and $d:M^2\rightarrow \Sigma_Z$ given by the \textbf{symmetric difference} of sets
\[
d(A,B):=(A\cup B)\backslash(A\cap B)=(A\cup B)\cap(A\cap B)^c=(A\cap B^c)\cup(B\cap A^c)=A\backslash B~\cup~B\backslash A.
\]
Let $A,B\in M$. If there exists $x\in Z$ such that $\{x\}\in\Sigma_Z$, then $B_{\{x\}}(A)=\{A\}$ (i.e., the $d$-topology is discrete since $A\in M$ was arbitrary). Similarly, if $\{x,y\}\in\Sigma_Z$, then
\[
B_{\{x,y\}}(A)\subset\left\{
                       \begin{array}{ll}
                         \big\{A,A\backslash\{x\},A\backslash\{y\}\big\}, & \textrm{if $x,y\in A$,} \\
                         \big\{A,A\backslash\{x\},A\cup\{y\}\big\}, & \textrm{if $x\in A$, $y\in A^c$,} \\
                         \big\{A,A\cup\{x\},A\cup\{y\}\big\}, & \textrm{if $x,y\in A^c$.}
                       \end{array}
                     \right.
\]
More generally, for any $\emptyset\neq \varepsilon\in\Sigma_Z$,
\[
B_\varepsilon(A)\subset\{(A\backslash\alpha)\cup \beta:\alpha\subset\varepsilon\cap A,\beta\subset\varepsilon\cap A^c,\alpha\cup\beta\subsetneq\varepsilon\}.
\]
\end{example}

\begin{dfn}[{sv-metrizable space}]\label{svMetrizDfn}
A space is \textbf{sv-metrizable} if it is homeomorphic to a sv-metric space.
\end{dfn}

A metric can be viewed as a sv-metric, and so it is clear that every metrizable space is sv-metrizable. Whether or not a given topological space is sv-metrizable (Question \ref{GHDQuestion3}) depends on the severity of the restrictions (in Definition \ref{PartAlgDfn}) placed on the allowed class of partial algebras. 

\subsection{Postmeasures and decomposition of $d_H$}
A treatment of ordered groups, rings, and modules (Definition \ref{OrdAbelGrp}) can be found in \cite{steinberg2010}.

\begin{dfn}[{Ordered abelian group, Ordered ring, Ordered module, Ordered algebra}]\label{OrdAbelGrp}
A \textbf{(partially) ordered abelian group} is a poset $A=(A,\leq)$ such that (i) $A$ is an abelian group and (ii) for all $a,b,c\in A$, if $a\leq b$, then $a+c\leq b+c$.

The subset $A^+=\{a\in A:a\geq 0\}$ is called the ``\textbf{partial order}'' of $(A,\leq)$, because it determines $(A,\leq)$ through the following properties: $A^+$ satisfies ``$0\in A^+$'' by reflexivity of $\leq$, ``$a,-a\in A^+$ $\Rightarrow$ $a=0$'' by antisymmetry of $\leq$, and ``$A^++A^+:=\{a+b:a,b\in A^+\}\subset A^+$'' by transitivity of $\leq$ (\textbf{footnote}\footnote{\textbf{How $A^+$ determines $(A,\leq)$:} Let $A$ be an abelian group. Fix a subset $P\subset A$ such that (i) $0\in P$, (ii) $a,-a\in P$ $\Rightarrow$ $a=0$, and (iii) $P+P:=\{a+b:a,b\in P\}\subset P$. For any $a,b\in A$, define $a\leq b$ by $b-a\in P$ (conventionally written as $a-b\geq 0$, so that $P=\{a\in A:a\geq 0\}$). Then, since $a-b=(a+c)-(b+c)$, we see that for all $a,b,c\in A$, $a\leq b$ $\iff$ $a+c\leq b+c$. Also, property (i) makes $\leq$ reflexive, property (ii) makes $\leq$ antisymmetric, and property (iii) makes $\leq$ transitive.
}).

An \textbf{ordered ring} is a poset $R=(R,\leq)$ such that (i) $R$ is a ring, (ii) $(R,+)=(R,+,\leq)$ is an ordered abelian group, and (iii) $R^+$ is closed under mutiplication (i.e., $R^+R^+:=\big\{\textrm{finite sums}~\sum_{i=1}^ka_ib_i:a_i,b_i\in R^+,k\in\mathbb{N}\big\}\subset R^+$). Given an ordered ring $R$, an \textbf{ordered $R$-module} (resp., \textbf{ordered $R$-algebra}) is a poset $M=(M,\leq)=(_RM,\leq)$ such that (i) $M$ is an $R$-module (resp., $R$-algebra) and (ii) $R^+M^+:=\big\{\textrm{finite sums}~\sum_{i=1}^ka_im_i:a_i\in R^+,m_i\in M^+,k\in\mathbb{N}\big\}\subset M^+$ (resp., and $(M,\leq)$ is an ordered ring). (\textbf{footnote}\footnote{In this definition, $R^+R^+$ and $R^+M^+$ (as defined) can be replaced with $R^+\ast R^+:=\{ab:a,b\in R^+\}$ and $R^+\ast M^+:=\{am:a\in R^+,m\in M^+\}$ respectively. This is because we already have $R^++R^+\subset R^+$ and $M^++M^+\subset M^+$ respectively.
})
\end{dfn}
Observe that if $a,b,c,d\in A$ are such that $a\leq b$ and $c\leq d$, then $a+c\leq b+c\leq b+d$. Also, if $a,b,c,d\in A$ are such that $a<b$ and $c\leq d$, then ($a+c\leq b+c\leq b+d$, where $a+c\neq b+c$, and so) $a+c<b+c\leq b+d$. If $a,b\in R$ are such that $a\leq b$, then (i) $-a\geq -b$ (as $c\geq 0$ iff $0\geq -c$) and (ii) for any $p\in R^+$ we have $ap\leq bp$, $pa\leq pb$ (because $a\leq b$ $\iff$ $0\leq -a+b$, $\Rightarrow$ $0\leq (-a+b)p=-ap+bp$ and $0\leq p(-a+b)=-pa+pb$). Similarly, if $m,n\in M={}_RM$ are such that $m\leq n$, then (i) $-m\geq -n$ and (ii) for any $p\in R^+$ we have $pm\leq pn$.

\begin{dfn}[{Postmeasure, Postmeasure space}]\label{PostMeasDfn}
Let $Z$ be a set, $(\Sigma_Z,\preceq,\uplus)$ a partial algebra on $Z$, and $K=(K,\leq,+)$ an ordered abelian group. A map $\mu:(\Sigma_Z,\preceq,\uplus)\rightarrow (K,\leq,+)$ is a \textbf{postmeasure} on $Z$ (making $~Z=(Z,\Sigma_Z,\mu)~$ a \textbf{postmeasure space}) if:
\begin{enumerate}
\item $\mu(R)=0$ $\iff$ $R=0$. (faithfulness)
\item $R\prec S$ implies $\mu(R)\leq \mu(S)$. (monotonicity)
\item $\mu(R\uplus S)\leq \mu(R)+\mu(S)$. (subadditivity)
\end{enumerate}
\end{dfn}

\begin{thm}\label{dHausCompThm}
Let $M$ be a set. (i) If $d_{sv}:M^2\rightarrow\big(\Sigma_Z,\preceq,\uplus\big)$ is a sv-metric and $\mu:(\Sigma_Z,\preceq,\uplus)\rightarrow (\mathbb{R},\leq,+)$ a real-valued postmeasure, then the composition $\mu\circ d_{sv}:M^2\rightarrow\mathbb{R}$ is a metric.

(ii) If $\mu$ is strictly monotone {\small (i.e., $R\prec S$ implies $\mu(R)<\mu(S)$) and $\mu^{-1}(r)\backslash\{0\}\neq\emptyset$ for all $r>0$}, then $\tau_{\mu\circ d_{sv}}\subset\tau_{d_{sv}}$ (hence $\tau_{d_{sv}}$ is Hausdorff).

(iii) There exists a postmeasure $\mu$ and a sv-metric $d_{sv}$ such that
\begin{align}
d_H=\mu\circ d_{sv},~~~~\textrm{(i.e., $d_H$ is in the postmeasure family, Definition \ref{PostMeasFam}).}
\end{align}
 \end{thm}
\begin{proof}
(i) This is immediate by construction (involving a simple verification of the axioms in Definition \ref{MetricDfn}). (ii) If $\mu$ is strictly monotone and $\mu^{-1}(r)\backslash\{0\}\neq\emptyset$ for all $r>0$, then
\[
0\prec\varepsilon\in\mu^{-1}(r)~\textrm{~}~\Rightarrow~\textrm{~}~B_\varepsilon^{d_{sv}}(m)\subset B_r^{\mu\circ d_{sv}}(m).
\]
(iii) By (i), we only need to specify the maps $\mu$ and $d_{sv}$.
Let $X$ be a metric space and $M:=BCl(X)\subset\mathcal{P}^\ast(X)$. Then, with $Z:=\mathbb{R}_+:=[0,\infty)$, $\Sigma_Z:=\mathcal{P}^\ast(Z)$, ~$\preceq\!~:=~\!\leq_{\sup}~$ (in the sense that ``$R\leq_{\sup} S$ $\iff$ $\sup(R)\leq\sup(S)$''), ~$\uplus:=+_{\textrm{cpwise}}$~ as componentwise addition $R+_{\textrm{cpwise}}S:=\{r+s:r\in R,s\in S\}$, and $\mu:=\sup$, we have
\begin{align*}
d_H=\sup\circ d_{sv}&:BCl(X)^2\stackrel{d_{sv}}{\longrightarrow}\Big(\mathcal{P}^\ast(\mathbb{R}_+),\leq_{\sup},+_{\textrm{cpwise}},\mathcal{P}^\ast(\mathbb{R}_+)\Big)\stackrel{\sup}{\longrightarrow}(\mathbb{R},\leq,+),\\
d_{sv}(A,B)&:=\{\textrm{dist}(a,B):a\in A\}\cup\{\textrm{dist}(b,A):b\in B\}~~\\
&=\{|\textrm{dist}(x,A)-\textrm{dist}(x,B)|:x\in A\cup B\}\\
&~~~~\textrm{or}~~\{|\textrm{dist}(x,A)-\textrm{dist}(x,B)|:x\in X\}.
\end{align*}
The triangle inequality holds in the form $~d_{sv}(A,B)\leq_{\sup}d_{sv}(A,C)+_{\textrm{cpwise}} d_{sv}(C,B)$.

\emph{Alternatively} (observing that ~$\textrm{dist}(a,B)=\inf\{r>0:a\in B_r\}=\sup\big(\{r>0:a\in B_r\}^c\big)$,~ where $A_r=\overline{N}_r(A):=\{x\in X:\textrm{dist}(x,A)\leq r\}$), we have
\begin{align*}
&d_H=\sup\circ d^{sv}: BCl(X)^2\stackrel{d^{sv}}{\longrightarrow}\Big(\mathcal{P}^\ast(\mathbb{R}_+),\leq_{\sup},+_{\textrm{cpwise}},\mathcal{P}^\ast(\mathbb{R}_+)\Big)\stackrel{\sup}{\longrightarrow}(\mathbb{R},\leq,+),
\end{align*}
where
\[
d^{sv}(A,B):=\{r>0:A\cup B\subset A_r\cap B_r\}^c=\{r>0:A\subset B_r\}^c\cup \{r>0:B\subset A_r\}^c,
\]
with the triangle inequality again in the form $~d^{sv}(A,B)\leq_{\sup}d^{sv}(A,C)+_{\textrm{cpwise}} d^{sv}(C,B)$.
\end{proof}

\begin{dfn}[{Postmeasure family of distances}]\label{PostMeasFam}
Let $X$ be a space. A distance $d:Cl(X)^2\rightarrow\mathbb{R}$ belongs to the \textbf{postmeasure family} if it factors (as in Theorem \ref{dHausCompThm}) in the form $d=\mu\circ d_{sv}$ for a postmeasure $\mu$ and a sv-distance $d_{sv}$.
\end{dfn}

In addition to the sv-metric encountered in the proof of Theorem \ref{dHausCompThm}(iii), we will explicitly describe another sv-metric in Definition \ref{PractSvMetEx1}. More examples of sv-metrics can be readily obtained by factoring the generalized Hausdorff distances $d_{p,\nu}$, $d_{p,\nu,\mu}$, $d_{p,\nu,\mu,\alpha,\beta}$, $d_{F,G}$ of Examples \ref{AvgHausDist}, \ref{DIntHausDist}, \ref{MetVecSpDist}, \ref{FunHausDist} respectively (as done for $d_H$ in Theorem \ref{dHausCompThm}(iii)).

\begin{example}[{Multipath-valued metric}]\label{PractSvMetEx1}
Let $X$ be a metric space and $Y:=[0,1]^n\subset\mathbb{R}^n$. Let two $n$-dimensional objects $A,B\in BCl_Y(X)$ be given (in parameterized form) by the closures of images of bounded maps $f,g:Y\rightarrow X$, that is, $A=q(f)$ and $B=q(g)$ as elements of the metric space $\big(BCl_Y(X),d_H\big)$. Given positive integers $n$ and $k$, let $m:=k^n$,
\[
Y_m:=\{0=0/k,1/k,2/k,...,k/k=1\}^n\subset{Y},
\]
and consider the discrete approximations of $A$ and $B$ by the finite sets $A_m:=g({Y_m})$ and $B_m:=g({Y_m})$ in
\[
FS_m(X):=BCl_{Y_m}(X)=\{C\in BCl_Y(X):1\leq|C|\leq m\}.
\]

We note that the discrete approximation here is only introduced for easy visualization and computation, otherwise the approximation is not strictly essential in the discussion. With the set
\[
Z:=\textrm{Paths}(FS_m(X)):=C([0,1],FS_m(X))
\]
of paths in $(FS_m(X),d_H)$, and a suitable subset $\Sigma_Z\subset\mathcal{P}^\ast(Z)$, define an sv-metric by
\begin{align*}
d_{sv}^m:BCl_Y(X)^2&\rightarrow \Sigma_Z\subset\Big(\mathcal{P}^\ast(Z),\preceq,\uplus\Big),~(A,B)\mapsto~\textrm{Paths}_{A_m,B_m}(FS_m(X)),
\end{align*}
where $\textrm{Paths}_{A_m,B_m}(FS_m(X))$ is taken to be a subcollection of paths in $(FS_m(X),d_H)$ from $A_m$ to $B_m$, the ordering $\preceq$ is given by
\[
\Gamma\preceq \Gamma'~~\iff~~l_{\min}(\Gamma)\leq l_{\min}(\Gamma'),~~\textrm{~}~~l_{\min}(\Gamma):=\inf\big\{\textrm{length}(\gamma):\gamma\in \Gamma\big\},
\]
and the operation $\uplus:=\textrm{conc}_{\textrm{ptwise}}$ is ``pointwise'' path concatenation, in the sense that
\[
\Gamma\uplus\Gamma':=\{\gamma\cdot\gamma':\gamma\in\Gamma,\gamma\in\Gamma'\},~~\textrm{where}~~
(\gamma\cdot\gamma')(t):=
\left\{
  \begin{array}{ll}
    \gamma(2t), & t\in[0,1/2], \\
    \gamma'(2t-1), & t\in[1/2,1].
  \end{array}
\right.
\]
This concludes Example \ref{PractSvMetEx1}.
\end{example}

\begin{rmk}\label{PracSvMetRmk}
From Example \ref{PractSvMetEx1}, we obtain (by  Theorem \ref{dHausCompThm}) the associated real-valued distance
\[
d_m=l_{\min}\circ d_{sv}^m: BCl_Y(X)^2\stackrel{d_{sv}^m}{\longrightarrow}\Big(\Sigma_{Z},l_{\min},\textrm{conc}_{\textrm{ptwise}}\Big)\stackrel{l_{\min}}{\longrightarrow}(\mathbb{R},\leq,+),
\]
where the additivity of path length over path concatenation (i.e., $\textrm{length}(\gamma\cdot\eta)=\textrm{length}(\gamma)+\textrm{length}(\eta)$ ) ensures that ``minimum path length'' $l_{\min}$ is a postmeasure. Therefore, whenever $\big(FS_m(X),d_H)$ is a \emph{length space}, in the sense that $d_H(A_m,B_m)=l_{\min}\big(\textrm{Paths}_{A_m,B_m}(FS_m(X))\big)$, we get
\[
d_m(A,B)=d_H(A_m,B_m).
\]

This example therefore points to another way of decomposing $d_H$ that is different from that encountered in Theorem \ref{dHausCompThm}(iii). In general, the metrics $d_m$ are approximations (or ``generalizations'') of $d_H$.
\end{rmk}

\section{\textnormal{\bf Generalized Hausdorff distances (ghd's)}} \label{GenHausSec}
\noindent We will describe two main classes of generalized Hausdorff distances, namely, the \emph{relational} and \emph{integral} classes, which overlap since Hausdorff distance can be recovered from each (Remarks \ref{HDRremark2}, \ref{HDRremark3}, \ref{HDRremark4}).

\begin{dfn}[{Class of generalized Hausdorff distances (ghd-class)}]\label{HausDistClass}
Let $X$ be a space. A family of distances $\{d_\lambda:Cl(X)^2\rightarrow\mathbb{R}\}_{\lambda\in\Lambda}$ is a \textbf{ghd-class} on $Cl(X)$ if the following hold:
\begin{enumerate}
\item If $X$ is a metric space, then there exists $\lambda\in\Lambda$ such that $d_\lambda=d_H$.
\item There is a fixed $Cl(X)$-dependent topology on $\Lambda$ with respect to which the map $d:\Lambda\rightarrow C(Cl(X)^2,\mathbb{R}),\lambda\mapsto d_\lambda$ is continuous (where $C(U,V)$ denotes a space of continuous functions from $U$ to $V$). ({\bf footnote}\footnote{
    Here, continuity of $\lambda\mapsto d_\lambda$ justifies the use of ``generalized'' in ``generalized Hausdorff distances''.
    })
\end{enumerate}
\end{dfn}

\subsection{Relational ghd-classes} \label{RelHausSec}

A ghd-class $\{d_\lambda\}_{\lambda\in\Lambda}$, on $Cl(X)$, (Definition \ref{HausDistClass}) is
\textbf{relational} (or \textbf{relation-based}) if for any $\lambda\in \Lambda$ and any $A,B\in Cl(X)$, $d_\lambda(A,B)$ is defined using a relation $R_\lambda\subset A\times B$. Recall from Remark \ref{SemiMetConv} that semimetrics can be converted to pseudometrics, which in turn can be converted to metrics by ``quotienting''. Some of the distances below will be semimetrics or pseudometrics.

\begin{example}[{Upper relational (UR-) Hausdorff distances}]\label{RelHausDist1}
Let $X$ be a metric space, $A,B\subset X$, and $R\subset A\times B$ (a \textbf{relation} between $A$ and $B$). The \textbf{domain} and \textbf{range} of $R$ are respectively
\begin{align*}
\textrm{dom}(R)&=A_R:=\{a\in A:\textrm{some}~(a,b)\in R\}~~\textrm{~}~~\textrm{and}~~\\
\textrm{ran}(R)&={}_RB:=\{b\in B:\textrm{some}~(a,b)\in R\}.
\end{align*}
We say a relation $R\subset A\times B$ is \textbf{complete} if $A_R=A$ and ${}_RB=B$, and $R$ is \textbf{intersection-complete} if $R$ is the intersection of some family $\{R_\alpha\}_{\alpha\in\mathcal{A}}$ of complete relations $R_\alpha\subset A\times B$. Let $T(X):=BCl(X)^2$ and
{\small
\[
\textstyle\textrm{ICR}(X):=\bigcup_{A,B\in BCl(X)}\{\textrm{intersection-complete relations in~}~\mathcal{P}(A\times B)\},
\]}where $\mathcal{P}(A\times B)=\bigcup_{C\in\mathcal{P}(A)}\mathcal{P}(B)^C$, with each relation $R\subset A\times B$ viewed as the map
\[
f_R:\textrm{dom}(R)\subset A\rightarrow \mathcal{P}(B),~a\mapsto \{b\in B:(a,b)\in R\}.
\]
Consider any ``\textbf{selection}'' {\small $R\in \textrm{ICR}(X)^{T(X)}$} of intersection-complete relations (i.e., {\small $R:T(X)\rightarrow \textrm{ICR}(X),~(A,B)\mapsto R(A,B)\subset A\times B$}) satisfying {\small $R(A,A)=\{(a,a):a\in A\}$} for all $A$. (\textbf{footnote}\footnote{For example, $R(A,B)\subset A\times B$ could be the graph $G(f^{R(A,B)})$ of a \emph{partial function} $f^{R(A,B)}:\textrm{dom}(R(A,B))\rightarrow B$.}). To $R$ we associate a pseudometric (call it the \textbf{upper relational (UR$(R)$-) Hausdorff distance}) on $BCl(X)$ given by
\begin{align}
\label{RelHausDistEq1}\textstyle d^R(A,B):=\sup\{d(a,b):(a,b)\in R(A,B)\},
\end{align}
where the \textbf{triangle inequality (TI)} for $d^R$ requires the following \emph{necessary and sufficient} condition on $R$.
\end{example}

\begin{prp}[{TI-criterion}]\label{UrbGhdTIC1}
The triangle inequality holds for $d^R$ in (\ref{RelHausDistEq1}) if and only if
 $R$ satisfies the following criterion for any $A,B,C\in BCl(X)$.

\textbf{TI-criterion}: For any $(a,b)\in R(A,B)$ and $\varepsilon>0$, there exist (possibly $(a,b)$-dependent) $(a_\varepsilon,c_{\varepsilon 1})\in R(A,C)$ and $(c_{\varepsilon 2},b_\varepsilon)\in R(C,B)$ such that $d(a,b)\leq d(a_\varepsilon,c_{\varepsilon 1})+d(c_{\varepsilon 2},b_\varepsilon)+\varepsilon$. (\textbf{footnote}\footnote{
For the latter inequality, it is sufficient, but perhaps not necessary, to require that ``$d(a,a_\varepsilon)<\varepsilon/3$, $d(c_{\varepsilon 1},c_{\varepsilon 2})<\varepsilon/3$, and $d(b_\varepsilon,b)<\varepsilon/3$'', since the triangle inequality for $X$ implies $d(a,b)\leq d(a,a_\varepsilon)+d(a_\varepsilon,c_{\varepsilon 1})+d(c_{\varepsilon 1},c_{\varepsilon 2})+d(c_{\varepsilon 2},b_\varepsilon)+d(b_\varepsilon,b)$.
}).
\end{prp}
\begin{proof}
Suppose $d^R(A,B)\leq d^R(A,C)+d^R(C,B)$. Fix $(a,b)\in R(A,B)$ and $\varepsilon>0$. Then there exist $(a_\varepsilon,c_{\varepsilon1})\in R(A,C)$ and $(c_{\varepsilon2},b_\varepsilon)\in R(C,B)$ such that $d^R(A,C)<d(a_\varepsilon,c_{\varepsilon1})+\varepsilon/2$ and $d^R(C,B)<d(c_{\varepsilon2},b_\varepsilon)+\varepsilon/2$, and so
$d(a,b)\leq d^R(A,B)<d(a_\varepsilon,c_{\varepsilon1})+d(c_{\varepsilon2},b_\varepsilon)+\varepsilon$.

Conversely, suppose the TI-criterion holds. Fix $(a,b)\in R(A,B)$ and $\varepsilon>0$. Then there exist $(a_\varepsilon,c_{\varepsilon1})\in R(A,C)$ and $(c_{\varepsilon2},b_\varepsilon)\in R(C,B)$ such that
$d(a,b)<d(a_\varepsilon,c_{\varepsilon1})+d(c_{\varepsilon2},b_\varepsilon)+\varepsilon$. So, $d^R(A,B)\leq d(a,b)<d(a_\varepsilon,c_{\varepsilon1})+d(c_{\varepsilon2},b_\varepsilon)+\varepsilon\leq d^R(A,C)+d^R(C,B)+\varepsilon$.
\end{proof}

\begin{rmk}[{Recovering $d_H$ from $d^R$}]\label{HDRremark2}
If we define $R_H\in \textrm{ICR}(X)^{T(X)}$ by
\begin{align*}
R_H(A,B):&\textstyle =\{(a,b)\in A\times B:d(a,b)\leq d_H(A,B)\}\\
&\textstyle=\bigcap_{\lambda>1}\{(a,b)\in A\times B:d(a,b)\leq \lambda d_H(A,B)\},
\end{align*}
then we get
\[
d^{R_H}(A,B)=d_H(A,B).
\]
\end{rmk}

\begin{example}[{Collective upper relational (CUR-) Hausdorff distances}]\label{RelHausDist2}
Let $\mathcal{R}\subset \textrm{ICR}(X)^{T(X)}$, define
\begin{align}
\label{RelHausDistEq2}\textstyle d^{\mathcal{R}}(A,B):=\inf\{d^R(A,B):R\in\mathcal{R}\}=\inf_{R\in\mathcal{R}}\sup_{(a,b)\in R}d(a,b),
\end{align}
and call $d^\mathcal{R}$ the \textbf{collective upper relational (CUR$(\mathcal{R})$-) Hausdorff distance} on $BCl(X)$.
\end{example}

\begin{prp}\label{UrbGhdTIC2}
If the $d^R$'s satisfy the triangle inequality, then so does $d^\mathcal{R}$.
\end{prp}
\begin{proof}
Fix any $\varepsilon>0$. Choose $R_\varepsilon,R'_\varepsilon,R''_\varepsilon\in\mathcal{R}$ such that $R_\varepsilon\subset R'_\varepsilon\cap R''_\varepsilon:=\{R'_\varepsilon(A,B)\cap R''_\varepsilon(A,B):A,B\in BCl(X)\}$ [ where $R\subset S$ iff $R(A,B)\subset S(A,B)$ for all $A,B\in BCl(X)$ ], $d^{R'_\varepsilon}(A,C)<d^{\mathcal{R}}(A,C)+\varepsilon/2$, and $d^{R''_\varepsilon}(C,B)<d^{\mathcal{R}}(C,B)+\varepsilon/2$. Then
{\small $d^{\mathcal{R}}(A,B)\leq d^{R_\varepsilon}(A,B)\leq d^{R_\varepsilon}(A,C)+d^{R_\varepsilon}(C,B)\leq d^{R'_\varepsilon}(A,C)+d^{R''_\varepsilon}(C,B)<d^{\mathcal{R}}(A,C)+d^{\mathcal{R}}(C,B)+\varepsilon$}.
\end{proof}

\begin{example}[{Lower relational (LR-) Hausdorff distances}]\label{RelHausDist3}
Consider the setup in Example \ref{RelHausDist1}. With $R(X):=\{\textrm{relations}~R\subset A\times B,~~A,B\in BCl(X)\}$, let $S\in R(X)^{T(X)}$ be a selection of relations such that ``$S(A,A)$ is symmetric for all $A$'' (or more strongly if need be, ``$S(A,B)^{-1}=S(B,A)$ for all $A,B$''). Fix $A,B\in BCl(X)$. Let $A_S=A_{S(A,B)}$ and ${}_SB={}_{S(A,B)}B$ be the domain and range of $S(A,B)$ as before. Define the semimetric (call it the \textbf{lower relational (LR$(S)$-) Hausdorff semimetric}) on $BCl(X)$ given by
\begin{align}
\label{RelHausDistEq3}d_S(A,B)&\textstyle:=d_H(A_S,{}_SB)=\inf\{r>0:A_S\subset\overline{N}_r({}_SB),{}_SB\subset\overline{N}_r(A_S)\}\nonumber\\
&\textstyle=\max\{\sup_{b\in\!\!~{}_SB}\textrm{dist}(A_S,b),\sup_{a\in A_S}\textrm{dist}(a,{}_SB)\}\nonumber\\
&\textstyle=\sup_{x\in A_S\cup{}_S\!B}|\textrm{dist}(x,A_S)-\textrm{dist}(x,{}_SB)|\nonumber\\
&\textstyle =\sup_{x\in X}|\textrm{dist}(x,A_S)-\textrm{dist}(x,{}_SB)|.
\end{align}

\begin{prp}\label{UrbGhdTIC3}
The distance $d_S$ satisfies the triangle inequality if we choose $S\in R(X)^{T(X)}$ such that ~${}_{S(A,B)}B=B_{S(B,C)}$,~ for all ~$A,B,C\in BCl(X)$.
\end{prp}

\begin{rmk}[{Recovering $d_H$ from $d_S$}] \label{HDRremark3}
If we define $S_H\in R(X)^{T(X)}$ to be such that $S_H(A,B)$ is complete for each pair $(A,B)$, then
\[
d_{S_H}(A,B)=d_H(A,B).
\]
\end{rmk}
\end{example}

\begin{rmk}[Adaptability]\label{AdaptRmk}
Depending on specific application requirements, the supremum and/or infimum in Eqn. (\ref{RelHausDistEq1}), (\ref{RelHausDistEq2}), or (\ref{RelHausDistEq3}) may be replaced with more suitable operations, which include measures and integrals. Moreover, each class of subsets $\mathcal{S}(X)\subset BCl(X)$ can admit distances that are more suitable/convenient for that class than for all of $BCl(X)$. These remarks also accordingly apply to the integral ghd's of Section \ref{IntHausSec}.
\end{rmk}

\subsection{Integral ghd-classes} \label{IntHausSec}
A ghd-class $\{d_\lambda\}_{\lambda\in\Lambda}$, on $Cl_Y(X)$, (Definition \ref{HausDistClass}) is
\textbf{integral} (or \textbf{integration-based}) if for each $\lambda\in\Lambda$, $d_\lambda$ is defined using a real-valued integration functional $\int_{X\times Y_\lambda^2}(~)d\mu_\lambda$ involving a function space $\mathcal{F}_\lambda\subset X^Y$ (for some set $Y$), where $Cl_Y(X)$ is a quotient space of $\mathcal{F}_\lambda$ in the sense of \cite[Theorem 5.1]{akofor2025}.

In this section, the operation $q(\cdot)$ denotes the mapping in Equation (\ref{UnordMapEq}). In all expressions that follow, we assume the functions involved are sufficiently measurable/continuous if need be, and that the measures are appropriately chosen to ensure that the relevant metric axioms in Definition \ref{MetricDfn} hold.

\begin{rmk}
Let $X$ be a space and $Y$ a set. Note that a metric $d(q(f),q(g))$ on $Cl_Y(X)$ corresponds to the pseudometric $d'(f,g):=d(q(f),q(g))$ on $X^Y$.
\end{rmk}

Let $X$ be a metric space, $Y$ a set, and $f,g\in BCl(Y,X):=q^{-1}(BCl_Y(X))$. Recall from Equation (\ref{HausDistDef}) that the Hausdorff distance between $q(f)$ and $q(g)$ in $BCl_Y(X)$ takes the form
\[
\textstyle d_H\left(q(f),q(g)\right)=\max\big\{\sup_y\inf_{y'}d(f(y),g(y')),\sup_{y'}\inf_yd(f(y),g(y'))\big\}.
\]
Consider measure spaces $(X,\nu)$ and $(Y,\mu)$. For some positive function $\rho:X\rightarrow(0,\infty)$, let (\textbf{footnote}\footnote{
Here, the map $\chi_A:X\rightarrow\mathbb{R}$ given by $\chi_A|_{A}=1$ and $\chi_A|_{X\backslash A}=0$ is the \textbf{characteristic function} of the set $A\subset X$.
})
\[
{\rho}_{fg}(x)\textstyle:=\chi_{q(f)\cup q(g)}(x){\rho}(x)~~\textrm{and}~~ \Lambda_{fg}(p):=
    \left(\int_{q(f)\cup q(g)}{\rho}(x)d\nu_x\right)^{1/p},~~ 1\leq p<\infty.
\]
Also recall the \textbf{uniform metric} on $BCl(Y,X)$ which is given by
\[
\textstyle d_u(f,g):=\sup_{y\in Y}d(f(y),g(y)).
\]

\begin{rmk}[{Recovering $d_H$ from integral ghd's}] \label{HDRremark4}
It is clear that $d_H$ can be recovered from the following distances.
\end{rmk}

\begin{example}[{$L^p(X,\nu)$-Integral Hausdorff distance}]\label{AvgHausDist}
We have the metrics:
{\small
\begin{align}
\label{PracGHDistEq1} d_{p,\nu}\left(q(f),q(g)\right)&\textstyle:=\left(\int_X|\textrm{dist}(x,q(f))-\textrm{dist}(x,q(g))|^p\rho_{fg}(x)d\nu_x\right)^{1/p},~\textrm{~}~1\leq p<\infty,\\
&\textstyle\leq \left(\int_X\sup_{y\in Y}|d(x,f(y))-d(x,g(y))|^p\rho_{fg}(x)d\nu_x\right)^{1/p}\nonumber\\
&\textstyle\leq \Lambda_{fg}(p)\!~d_u(f,g);\nonumber\\
d_{\infty,\nu}\big(q(f),q(g)\big)&:=\sup_{x\in X}|\textrm{dist}(x,q(f))-\textrm{dist}(x,q(g))|=d_H(q(f),q(g))\leq d_u(f,g).\nonumber
\end{align}}
\end{example}

\begin{example}[{$L^p(X\times Y^2,\nu\times\mu)$-Integral Hausdorff distance}]\label{DIntHausDist}
For sets $A,B\in Cl_Y(X)$, let {\small $\Delta(A,B):=(A\cup B)\backslash(A\cap B)=A\backslash B~\cup~B\backslash A$}. Consider real-valued functions $c_{fg}(x,y,y')$ on $X\times Y^2$ satisfying:
\begin{align}
\label{DIntHausEq1}
\left\{
  \begin{array}{ll}
    (i)~\textstyle  c_{fg}(x,y,y')=\chi_{\Delta\left(q(f),q(g)\right)}(x)c_{fg}(x,y,y'),~~\textrm{~}~~\forall x,y,y'. \\
 (ii)~\textstyle c_{fg}(x,y,y')>0,~\textrm{~}~\textrm{if $\Delta\left(q(f),q(g)\right)\neq\emptyset$ and $x\in \Delta\left(q(f),q(g)\right)$},~\textrm{~}~\forall y,y'.
  \end{array}
\right.
\end{align}
Suppose further that the following conditions are also satisfied:
\begin{align}
\label{DIntHausEq2}
\left\{
  \begin{array}{ll}
    (iii) & \int_{y}c_{fg}(x,y,y')\leq 1,~\forall x,y', \\
    (iv) & \int_{y'}c_{fg}(x,y,y')\leq1,~\forall x,y, \\
    (v) & \int_{y''}c_{fh}(x,y,y'')c_{hg}(x,y'',y')\geq c_{fg}(x,y,y'),~\forall x,y,y',
  \end{array}
\right.
\end{align}
where $\int_x(\cdots):=\int_X(\cdots) d\nu_x$ and $\int_y(\cdots):=\int_Y(\cdots) d\mu_y$. Let
\begin{align}
\label{PracGHDistEq2} d_{p,\nu,\mu}\left(q(f),q(g)\right)&\textstyle:=\left(\int_{X\times Y^2}|d(x,f(y))-d(x,g(y'))|^pc_{fg}(x,y,y')d\nu_xd\mu_{y,y'}\right)^{1/p}.
\end{align}
\end{example}

\begin{prp}\label{UrbGhdTIC4}
The distance in (\ref{PracGHDistEq2}) satisfies the triangle inequality (and so gives a metric). Moreover, the relations (\ref{DIntHausEq2}) can be modified such that, for a constant $K=K_{p,\nu,\mu}$,
\[
d_{p,\nu,\mu}\left(q(f),q(g)\right)\leq K\big[d_{p,\nu,\mu}\left(q(f),q(h)\right)+d_{p,\nu,\mu}\left(q(h),q(g)\right)\big].
\]
\end{prp}
\begin{proof}
If ~$d_{fg}(x,y,y'):=|d(x,f(y))-d(x,g(y'))|$, then
\begin{align*}
\textstyle\int_{x,y,y'}&\textstyle d_{fg}(x,y,y')c_{fg}(x,y,y')\stackrel{(\ref{DIntHausEq2})(v)}{\leq}\int_{x,y,y',y''}d_{fg}(x,y,y')c_{fh}(x,y,y'')c_{hg}(x,y'',y')\\
&\textstyle\leq\int_{x,y,y',y''}[d_{fh}(x,y,y'')+d_{hg}(x,y'',y')]c_{fh}(x,y,y'')c_{hg}(x,y'',y')\\
&\textstyle=\int_{x,y,y',y''}d_{fh}(x,y,y'')c_{fh}(x,y,y'')c_{hg}(x,y'',y')\\
&\textstyle~~+\int_{x,y,y',y''}d_{hg}(x,y'',y')c_{fh}(x,y,y'')c_{hg}(x,y'',y')\\
&\textstyle\stackrel{(\ref{DIntHausEq2})(iii),(iv)}{\leq}\int_{x,y,y''}d_{fh}(x,y,y'')c_{fh}(x,y,y'')+\int_{x,y',y''}d_{hg}(x,y'',y')c_{hg}(x,y'',y').
\end{align*}
\end{proof}

\begin{example}[{Weighted $L^p(X\times Y^2,\nu\times\mu)$-Integral Hausdorff distance}]\label{MetVecSpDist}
Let $X$ be a metric vector space and $p\geq 1$. Consider the same setup as in Example \ref{DIntHausDist} (including the functions $c_{fg}(x,y,y')$ and their properties (\ref{DIntHausEq1}) and (\ref{DIntHausEq2}) ). Pick positive real-valued functions $\alpha,\beta:X\times Y^2\rightarrow(0,\infty)$ and define (\textbf{footnote}\footnote{
Recall that given metrics $d_1,...,d_n:Z^2\rightarrow \mathbb{R}$, we get other metrics $d_p:Z^2\rightarrow\mathbb{R}$, $1\leq p\leq \infty$, given by
\[
\textstyle d_p(z,z'):=\left[\sum_id_i(z,z')^p\right]^{1/p}~\textrm{(for $1\leq p<\infty$)}~\textrm{and}~~d_\infty(z,z'):=\max_id_i(z,z').
\]
})
\begin{align}
\label{PracGHDistEq3} d_{p,\nu,\mu,\alpha,\beta}\left(q(f),q(g)\right)&\textstyle:=\Big(\int_{X\times Y^2}\big[\alpha(x,y,y')|d(x,f(y))-d(x,g(y'))|^p\\
&\textstyle~~+\beta(x,y,y') d\big(x-f(y),x-g(y')\big)^p\big]c_{fg}(x,y,y')d\nu_xd\mu_{y,y'}\Big)^{1/p}.\nonumber
\end{align}

Notice from this example that, on a metric space $X$, the presence of other structures (e.g., addition and scalar-multiplication) besides the metric can strictly increase the number of valid choices of distance between subsets of $X$.
\end{example}

\begin{example}[{Extended integral Hausdorff distance}]\label{FunHausDist}
Once again, consider the setup in Definition \ref{DIntHausDist}. Define
\begin{align}
\label{PracGHDistEq4} d_{F,G}\left(q(f),q(g)\right):=&\textstyle F\Big[\int_{X\times Y^2}G\big[d(x,f(y)),d(x,g(y')),d(f(y),g(y'))\big]\nonumber\\
&c_{fg}(x,y,y')d\nu_xd\mu_{y,y'}\Big],
\end{align}
for functions $F:[0,\infty)\rightarrow[0,\infty)$ and $G:\mathbb{R}^3\rightarrow[0,\infty)$ such that the following hold:
\begin{enumerate}
\item $d_{F,G}$ is well-defined (in the sense it does not depend on the choice of functions $f,g\in X^Y$ used to label elements of sets $A,B\in BCl(X)$ ), that is, if $q(h_1)=q(f)$ and $q(h_2)=q(g)$, then $~d_{F,G}\left(q(h_1),q(h_2)\right)=d_{F,G}\left(q(f),q(g)\right)$.
\item $F(t)=0$ $\iff$ $t=0$, and $~F(t+t')\leq F(t)+F(t')$.
\item ``$G(r,0,t)=0$ $\Rightarrow$ $r=0$'', ``$G(0,s,t)=0$ $\Rightarrow$ $s=0$'', and ``$G(0,0,t)=0$ $\iff$ $t=0$''.
\item $G(r,s,t)=G(s,r,t)$.
\item If need be, $F$ and $G$ are sufficiently elementary functions such that $d_{F,G}$ satisfies a(n approximate) triangle inequality.
\end{enumerate}

In this example, $d_{F,G}$ is in general a semimetric.
\end{example}

\section{\textnormal{\bf Conclusion and questions}}\label{DiscQuest}
\noindent We have seen (Theorem \ref{dHausCompThm}(iii)) that Hausdorff distance $d_H$ factors in the form $d_H=\mu\circ d_{sv}:M^2\stackrel{d_{sv}}{\longrightarrow}(\Sigma_Z,\preceq,\uplus)\stackrel{\mu}{\longrightarrow}(\mathbb{R},\leq,+)$, where $Z$ is a set, $\Sigma_Z\subset\mathcal{P}^\ast(Z)$ a partial algebra, $\mu$ a postmeasure, and $d_{sv}$ a set-valued metric. Based on this observation, we have constructed generalizations of $d_H$ (Equations (\ref{RelHausDistEq1}),(\ref{RelHausDistEq2}),(\ref{RelHausDistEq3}),(\ref{PracGHDistEq1}),(\ref{PracGHDistEq2}),(\ref{PracGHDistEq3}),(\ref{PracGHDistEq4})). The study of these more general distances should improve our understanding and interpretation of $d_H$ and its applications.

We remark here that any generalization of $d_H$ also has implications for $d_H$-dependent \emph{geometric} (i.e., \emph{geometrically invariant}) distances such as those in Definition \ref{GHDextends} below (also see \cite[including references 13-24 therein]{memoli2012} and \cite[Definition 3.4]{Gromov1999}).

\begin{dfn}[{Geometric distance between metric spaces}]\label{GHDextends}
Let Met be a (sub)category of metric spaces and consider metric spaces $X,Y,Z\in \textrm{Met}$. Let
\[
\textstyle E^Z(X,Y):=\{\textrm{simultaneous embeddings}~(I,J)~|~I:X\hookrightarrow Z,~J:Y\hookrightarrow Z\}\subset Z^X\times Z^Y
\]
and $\mathcal{G}\subset E(X,Y):=\bigcup_{Z\in\textrm{Met}}E^Z(X,Y)\subset \bigcup_{Z\in\textrm{Met}}Z^X\times Z^Y$. With $d_H=d_H^Z$ denoting Hausdorff distance between subsets of $Z$, the \textbf{$(H,\mathcal{G})$-geometric distance} (within Met) between $X$ and $Y$ is
\begin{align}
\textstyle d_{H,\mathcal{G}}(X,Y):=\inf\left\{d_H^Z(I(X),J(Y)):(I,J)\in \mathcal{G}\cap E^Z(X,Y),~Z\in \textrm{Met}\right\},
\end{align}
which is the smallest Hausdorff distance between simultaneous embeddings of $X$ and $Y$ from $\mathcal{G}$.

If only isometric embeddings and compact metric spaces are considered, then the resulting geometric distance is the \textbf{Gromov-Hausdorff distance}.

To obtain generalized versions of these geometric distances, we replace $d_H=d_H^Z$ in $d_{H,\mathcal{G}}$ with:
\begin{enumerate}
\item[(i)] $d^R=d^{R,Z}$ (from Eqn. (\ref{RelHausDistEq1})) to get $d^R_{\mathcal{G}}$ (call it UR$(R,\mathcal{G})$-\textbf{geometric distance}), where
\[
d^R_{\mathcal{G}}(X,Y):=\inf\left\{d^{R,Z}(I(X),J(Y)):(I,J)\in \mathcal{G}\cap E^Z(X,Y),~Z\in \textrm{Met}\right\};
\]
\item[(ii)] $d^\mathcal{R}=d^{\mathcal{R},Z}$ (from Eqn. (\ref{RelHausDistEq2})) to get $d^\mathcal{R}_{\mathcal{G}}$ (call it as CUR$(\mathcal{R},\mathcal{G})$-\textbf{geometric distance}), where
\[
d^\mathcal{R}_{\mathcal{G}}(X,Y):=\inf\left\{d^{\mathcal{R},Z}(I(X),J(Y)):(I,J)\in \mathcal{G}\cap E^Z(X,Y),~Z\in \textrm{Met}\right\};
\]
\end{enumerate}
and so on.
\end{dfn}

The following are some interesting questions.

\begin{question}\label{GHDQuestion3}
(i) From Definition \ref{svMetrizDfn}, do there exist spaces (resp., sv-metrizable spaces) that are not sv-metrizable (resp., not metrizable)?
(ii) How are sv-metric topologies (Definition \ref{SvOpenSetDfn}) related to topologies induced by uniform structures?
\end{question}

\begin{question}\label{GHDQuestion7}
In Example \ref{RelHausDist1}, when are the metrics $d_H$ and $d^R$ equivalent, i.e., when does there exist $c>0$ such that $d^R(A,B)/c\leq d_H(A,B)\leq c\!\!~d^R(A,B)$ for all $A,B$? The same can be asked of $d^{\mathcal{R}}$ (from Example \ref{RelHausDist2}) or of $d_S$ (from Example \ref{RelHausDist3}) in place of $d^R$.
\end{question}

\begin{question}\label{GHDQuestion8}
In Examples \ref{RelHausDist1} and \ref{RelHausDist3}, how are $d^R$ and $d_S$ related?
\end{question}

\begin{question}\label{GHDQuestion9}
Are all the integral ghd's of Section \ref{IntHausSec} relational ghd's?
\end{question}

\section*{\textnormal{\bf Declarations}}
The author has no competing interests to declare that are relevant to the content of this article.

\vspace{0.2cm}
\hrule
\endgroup
\end{document}